\documentclass{elsart}
\usepackage{amssymb}
\usepackage{amsfonts}
\usepackage{amsmath}
\journal{Journal of Mathematical Analysis and Applications}
\begin{document}

\begin{frontmatter}

\title{A Characterization of Norm Compactness in the Bochner Space $L^p \left(G ; B\right)$
For an Arbitrary Locally Compact Group $G$}

\author{Josh Isralowitz} \\

\address{Department of Mathematical Sciences, New Jersey Institute of Technology,
Newark, NJ 07102-1982 U.S.A. \thanksref{address} }

\thanks[address]{Current Address: Mathematics Department, University at Buffalo, Buffalo, NY 14260-2900}

\ead{jbi2@buffalo.edu}

\begin{abstract}
\indent In this paper, we generalize a result  of N. Dinculeanu
which characterizes
 norm compactness in the Bochner space $L^p\left(G ; B\right)$ in terms of an approximate identity
 and translation operators,
 where $G$ is a locally compact abelian group and $B$ is a Banach space.
 Our characterization includes the case where $G$ is nonabelian,
 and we weaken the hypotheses on the approximate identity used, providing
new results even for the case $B = \mathbb{C}$ and $G =
\mathbb{R}^n.$

\end{abstract}
\end{frontmatter}

\section{Introduction}
\indent The classical Kolmogorov-Riesz-Tamarkin  theorem
characterizes norm compact subsets $\Gamma$ of $L^p
\left(\mathbb{R}^n\right)$ in terms of the uniform convergence (in
the $L^p$ norm) of families of convolution and translation
operators, where a translation operator $T^h $ is defined by
$\left(T^h f\right) \left(x\right) = f\left(x - h\right)$ and a
convolution operator $J$ is given by $\left(J f\right)\left(x\right)
= \int_{\mathbb{R}^n} j\left(y\right) f\left(x - y\right) \, dy$ for
an appropriate kernel $j$ (for a brief history of this theorem see
\cite{Book2} or \cite{Book8}). Many authors have generalized this
theorem to $L^p\left(X\right)$ where $X$ is more general than
$\mathbb{R}^n$ and for function spaces that are more general than
the classical Lebesgue spaces $L^p$. One such generalization is the
following result, due to N. Dinculeanu \cite[Sect. 5]{Book2} :

\begin{thm}
\indent Let $G$ be a locally compact abelian group with Haar measure
$\mu$ and $B$ be a Banach space.
  A bounded subset $\Gamma$ of the Bochner space $L^p\left(G ; B\right), 1 \leq
p < \infty, $ is relatively norm compact
  iff

1) for any $A \subseteq G$ with $\mu\left(A\right) < \infty$, the
set $ \{\int_A f\left(y\right) \, d\mu\left(y\right)
: f \in \Gamma\} \subseteq B$ is relatively norm compact in $B$,\\

2) for any $\epsilon > 0$, there exists a compact $S \subseteq G$
such that \[\forall f \in \Gamma,  \ \ \ \left\Vert f\right\Vert _{p
; G \backslash S} < \epsilon,\]

and one of the following is true: \\

3a) let $\mathfrak{V}$ be a basis of relatively compact
neighborhoods of the identity in $G,$ and for each $V\in
\mathfrak{V}$ let $u_{V}$ be a positive, bounded, symmetric function
which vanishes outside $V$ where $\left\Vert u_{V}\right\Vert
_{1;G}=1.$ Then $\lim_{V\in \mathfrak{V}}\left\Vert
u_{V}\ast f-f\right\Vert  _{p;G}=0,$\\

3b) for any $\epsilon > 0$, there exists a neighborhood $\Theta$ of
the identity such that
\[\forall f \in \Gamma,
\ \ \ h \in \Theta \Longrightarrow \left\Vert T^h f - f\right\Vert  _{p ; G} < \epsilon.\]\\
\end{thm}

\indent Here we let
\[\left\Vert f \right\Vert _{p ; A} = \left(\int_A \left\Vert
f\left(x\right)\right\Vert _{B}^p \, d\mu \left(x\right)
\right)^\frac{1}{p}\] where $\left\Vert \cdot\right\Vert  _B$ is the
Banach space norm.

\indent We generalize the above result in two ways. First, we
generalize the theorem to include the case where $G$ is nonabelian.
Second, we weaken the hypotheses in property
 3a).  More specifically, we show that property 3a) can be replaced by the following:
for any $\epsilon > 0$ and compact $S \subseteq G$ provided by
property 2), there exists
 $j  \in
L^{p'}\left(G\right)$ such that $\forall f \in \Gamma, \left\Vert j
\ast f - f\right\Vert _{p ; S} < \epsilon$, where $p'$ is the
conjugate exponent of $p$ and where $j$ has compact support if $G$
is nonabelian or $p = 1$.
 We note that our weakening of these hypotheses
 appears to be new even for the classical cases $B = \mathbb{C}$ and $G = \mathbb{R}^n$.

\section{Preliminary Results}

In this section we provide some preliminary results the will be
needed to prove our main theorem. We start by introducing some
relevant definitions and notation. In all of the rest of the paper,
$G$ will denote a locally compact group with left invariant Haar
measure $\mu$. For $1 \leq p < \infty$ and a Banach space $B,
L^p\left(G; B\right)$ will denote the
 space of $\mu$ measurable functions $f : G
\rightarrow B$ satisfying $\left \Vert f \right \Vert_B ^p \in
L^1\left(G\right),$ where $L^p\left(G\right)$ denotes the standard
Lebesgue space of complex valued functions. The integral $\int_A
f(x) \, d\mu\left(x\right)$ over a set $A$ will also be written as
$\int_A f(x) \, dx$. We refer the reader to ~\cite[Chap. 3]{Book4}
for further properties of the Bochner integral.

We let $\triangle $ denote the continuous homomorphism from $G$ into
$\mathbb{R}^+$ defined by $\triangle\left(x\right) \mu\left(A\right)
= \mu\left(A x^{-1}\right)$, which by the definition of the Bochner
integral implies that $\int_G f\left(yx\right) \, dy = \triangle
\left(x\right) \int_G f\left(y\right) \, dy$.

If $j \in L^{p'}(G)$ and $f \in L^p(G, B)$, then the convolution
product $\left(j \ast f\right)\left(x\right)$ is defined by $
=\int_G j\left(xy\right) f\left(y^{-1}\right) \, dy.$ The
convolution can be defined equivalently  as $= \int_G
j\left(y\right)f\left(y^{-1}x\right) \, dy$, since the left
invariance of $\mu$ implies that $\int_G h\left(x y\right) \, dy =
\int_G h\left(y\right) \, dy$ for any $x \in G$ and any Bochner
integrable function $h$. ~\cite[Sect. 3]{Book2}

If we require that $j$ have compact support when either $G$ is
nonabelian or $p =1$, then the continuity of $\triangle$ and
H\"{o}lder's inequality implies that $j*f$ gives us a well defined
function from $G$ to $B$. Moreover, Lemma $7$ and the fact that $G$
is locally compact will imply that $j*f$ is continuous on $G$ so
that  $j\ast f \in L^r _{\rm loc}\left(G; B\right)$ for any $r \geq
1$, where $L^r _{\rm loc}\left(G; B\right)$ denotes the space of
$\mu$ measurable functions $f : G \rightarrow B$ such that $\left
\Vert f \right \Vert_{r; K} < \infty $ for every compact $K
\subseteq G$.

The main component in the proof of our main result will be Lemma 7.
Before proving this lemma, however, we need an elementary result
regarding the integration of functions in $\Gamma$ over sets of
small measure.

\begin{defn}  A subset  $\Gamma \subseteq L^1 \left(G ; B\right)$ is called
\textbf{Uniformly  Integrable}  on $S \subseteq G$ if for any
$\epsilon > 0$, there exists $\delta > 0$ such that \begin{gather*}
\forall f \in \Gamma, \ \ \ \
 A \subseteq S \mathrm{\ and \ } \mu(A) < \delta \Longrightarrow \| f \|_{1 ; A} < \epsilon.
\end{gather*}
\end{defn}

\begin{lem}
Let $\Gamma \subseteq L^1\left(G ; B\right)$ be bounded and let $S
\subseteq G$ be compact.  For any $\epsilon > 0$, assume that there
exists  $j \in
 L^{\infty}\left(G\right)$ with compact support $V$
   such that \[\forall f \in \Gamma, \ \ \
   \left\Vert f - j\ast  f\right\Vert _{1 ; S} < \epsilon.\]
Then $\Gamma ^{-}$ is uniformly integrable on $S$ where $\Gamma
^{-1} = \{f^{-} : f \in \Gamma\}$ and $f^{-}$ denotes the function
defined by $f^{-} \left(x\right) = f\left(x^{-1} \right).$

\end{lem}
\begin{pf}
First we show that for fixed $j$, the set
 $j \ast {\Gamma} = \{ j \ast f : f \in \Gamma\}$ is uniformly  integrable on $G$.
Note that as $V$ is compact, there exists $m > 0$ where
$\triangle\left(y^{-1}\right) > m$ for all  $y \in V$ since
$\triangle$ has a (nonzero) minimum on $V^{-1}.$ Let $\left\Vert
f\right\Vert _{\infty ; V}$ denote the essential supremum
$\left\Vert f\right\Vert _B$ over $V$.

Thus, we have
\begin{eqnarray} \left\Vert j \ast f\right\Vert _{1 ; A} & \leq &
\int_A \left(\int_{V} |j\left(x y\right)| \left\Vert f\left(y^{-1}\right)\right\Vert _B \, dy \right) \, dx \nonumber \\
& \leq & \left(\left\Vert j\right\Vert _{\infty ; V}\right)
\mu\left(A\right)
\int_{V} \left\Vert f\left(y^{-1}\right)\right\Vert _B \, dy \nonumber \\
& = & \left(\left\Vert j\right\Vert _{\infty ; V}\right)
\mu\left(A\right) \int_{V}
\left(\frac{\triangle\left(y^{-1}\right)}{\triangle\left(y^{-1}\right)}\right)
 \left\Vert f\left(y^{-1}\right)\right\Vert _B \, dy \nonumber \\ & < &
 \frac{1}{m} \left(\left\Vert j\right\Vert _{\infty ; V}\right)   \mu\left(A\right)
 \int_V \triangle\left(y^{-1}\right) \left\Vert f\left(y^{-1}\right)\right\Vert _B \, dy \nonumber \\ &
 \leq  & \frac{1}{m} \left(\left\Vert j\right\Vert _{\infty ; V}\right)    \mu\left(A\right)\int_{G} \left\Vert f\left(y\right)\right\Vert _B \,
dy \nonumber  \\ & \leq & M' \mu\left(A\right) \nonumber
\end{eqnarray} where $M'$ does not depend on $f \in \Gamma$, and so
 $j \ast {\Gamma} $ is uniformly integrable on $G$.

 \indent Now we show that $\Gamma$ is uniformly
integrable on $G$. By hypothesis, for any $\epsilon > 0$, we can
choose $j \in L^{\infty}\left(G\right)$
 with compact support, satisfying \[\forall f \in \Gamma, \ \ \
\left\Vert f - j\ast f\right\Vert _{1; G} \leq \frac{\epsilon}{2}.\]

\noindent Since $\{j\ast f : f \in \Gamma\}$ is uniformly integrable
on $G$, choose $\delta > 0$ so that \[\mu\left(A\right) < \delta
\Longrightarrow \forall f \in \Gamma, \ \ \ \left\Vert j\ast
f\right\Vert _{1; A} < \frac{\epsilon}{2}, \] and so
\[\mu\left(A\right) <
 \delta \Longrightarrow \forall f \in \Gamma, \ \ \ \left\Vert f\right\Vert _{1 ; A} \leq
\left\Vert f - j\ast f\right\Vert _{1; A} + \left\Vert j\ast
f\right\Vert _{1; A} < \epsilon.\]

\indent Finally, let $S \subseteq G$ be compact and let $A \subseteq
S$. By compactness, there exist $m$ and $M$ such that $0 < m <
\triangle\left(y^{-1}\right) < M $ for all  $ y \in S$. Thus, if $f
\in \Gamma ,$ we have \begin{eqnarray}\int_{A} \left\Vert
f\left(y^{-1}\right)\right\Vert _B \, dy &<& \frac{1}{m} \int_{A}
\triangle\left(y^{-1}\right) \left\Vert
f\left(y^{-1}\right)\right\Vert _B \, dy \nonumber
\\ & \leq& \frac{1}{m} \int_{G} \triangle\left(y^{-1}\right) \left(
\chi_{A^{-1}} \left(y^{-1}\right) \left\Vert f\left(y^{-1}\right)\right\Vert _B \right) \, dy \nonumber\\
&=& \frac{1}{m} \int_{G} \left( \chi_{A^{-1}}\left(y\right)\right)
\left\Vert f\left(y\right)\right\Vert _B \ \, dy \nonumber\\
&=& \frac{1}{m}\int_{A^{-1}} \left\Vert f\left(y\right)\right\Vert
_B \, dy. \nonumber\end{eqnarray}

 However, \begin{eqnarray}\mu\left(A^{-1}\right) &=& \int_G \chi_{A^{-1}} \left(x\right) \, dx \nonumber\\ &=& \int_G \chi_{A^{-1}}\left(x^{-1}\right)
\triangle\left(x^{-1}\right) dx\nonumber\\ & <& M \int_{G} \chi_{A^{-1}}\left(x^{-1}\right) \, dx \nonumber\\
&=& M \mu\left(A\right)\nonumber.
\end{eqnarray}   Therefore, since $\Gamma$ is
uniformly integrable on $G$,   we have that $\Gamma^{-}$ is
uniformly integrable on $S$.
\end{pf}\qed

Before we state and prove Lemma 7, which is crucial to our main
result, we will need two preliminary results, the first of which is
a norm compactness criteria for general normed spaces and was proven
in ~\cite[Lemma 3]{Book3}.

\begin{lem}
Let $S$ be any set, $F$ be a normed space with norm $\left\Vert \
\cdot \ \right\Vert _F ,$ and $ \{f_\alpha \}$ be a generalized
sequence of functions with $f_\alpha : S \rightarrow F$. \\
If $f: S \rightarrow F$ is a function such that $ \lim_{\alpha}
\left\Vert f_\alpha\left(s\right) - f\left(s\right)\right\Vert _F =
0$ uniformly for $s \in S$, and if each $f_\alpha\left(S\right)$ is
relatively norm compact in $F$, then $f\left(S\right)$ is relatively
norm compact in $F$.
\end{lem}

\indent The second preliminary result is the following form of the
Ascoli-Arzela theorem, which is a special case of the general
theorem ~\cite[Sect. 47]{Book6}. First we introduce some relevant
notation and definitions.

Let $S$ be a subspace of $G$. We let $C_B(S)$ denote the Banach
space of all continuous $B$ valued functions on $S$ endowed with the
``sup" norm $\mathrm{sup}_{x \in S} \|f(x)\|_B$.

\begin{defn}
If $S$ is a subspace of $G$, then a subset $K \subseteq C_B(S)$  is
said to be \textbf{equi-continuous at x} if for any $\epsilon > 0$,
there exists a relative neighborhood $N_x \subseteq S$ of $x$ such
that
\begin{align*}x' \in N_x \Longrightarrow \forall f \in K, \hspace{5mm} \|f(x) - f(x')\|_B < \epsilon.\end{align*}
If $K$ is equi-continuous on all of $S$, we say $K$ is
\textbf{equi-continuous} on $S$
\end{defn}

 \begin{lem} Let $S$ be a compact Hausdorff space and $B$ a Banach space.
A subset $K \subseteq C_B\left(S\right)$ is relatively norm
 compact iff

 1) $K$ is equi-continuous on $S$,\\

 2) for every $s \in S$, the set $K\left(s\right) = \{f\left(s\right) : f \in K\}$ is relatively
 norm compact in $B$.
 \end{lem}

\begin{lem}
Let $\Gamma \subseteq L^p\left(G ; B\right)$ be bounded and let $j
\in  L^{p'}\left(G\right)$ where  $j$ has compact support if either
$G$ is nonabelian or $p = 1$. Moreover, assume that when $A
\subseteq G $ and $\ \mu \left(A\right) < \infty$, then $\{\int_A
f^- \left(y\right) \, dy : f \in \Gamma\} \subseteq B$ is relatively
norm compact in $B$.

\indent Suppose that $S \subseteq G$ is compact and contains the
identity. Then $\left(j\ast\Gamma\right)|_S $ is relatively norm
compact in $C_B\left(S\right)$ where $\left(j\ast\Gamma\right)|_S =
\{\left(j \ast f\right)|_S : f \in \Gamma \}$ and $\left(j \ast
f\right)|_S$ is the restriction of $j \ast f$ to $S$.\end{lem}

\begin{pf} We first remark that we will show
the existence of each $\{\int_A f^- \left(y\right) \, dy : f \in
\Gamma\}$. We now check that $\left(j\ast\Gamma\right)$ is
equi-continuous on $S.$ Assume that $p > 1$, $G$ is nonabelian, and
supp $j = V$ is compact . Let $SV$ denote the set $\{ sv : s \in S,
v \in V\}$ so that
 $V \subseteq SV$ since $S$ contains the identity.
 Since $\Gamma$ is bounded and $SV$ is compact, we can choose $M > 0$ such that
\[ \forall f \in \Gamma, \ \ \    \left(\int_{SV} \left\Vert f\left(y^{-1}\right)
\right\Vert _B ^p \,
 dy \right)^\frac{1}{p}  < M.
\]

 Then using H\"{o}lder's equality, we compute that  \begin{eqnarray}\left\Vert \left(j\ast f\right)\left(x\right) -
\left(j\ast f\right)\left(x'\right)\right\Vert _B  &\leq& \int_V
|j\left(xy\right) - j\left(x' y\right)| \left\Vert
f\left(y^{-1}\right) \right\Vert _B \, dy \nonumber \\& \leq
&\left(\int_V |j\left(xy\right) - j\left(x' y\right)|^{p'} \, dy
\right)^\frac{1}{p'} \left( \int_V \left\Vert f\left(y^{-1}\right)
\right\Vert _B ^p \,
 dy \right)^\frac{1}{p} \nonumber \\& \leq & M \left(\int_V |j\left(xy\right) - j\left(x' y\right)|^{p'} \, dy
 \right)^\frac{1}{p'}. \nonumber
 \end{eqnarray}
 This clearly shows equi-continuity, since $j$ is translation
 continuous in the $L^{p'}$ norm. If $G$ is abelian and $p > 1$, then unimodularity
implies that $j$ does not need compact support.\\

\indent If $p = 1$ and $G$ is not necessarily abelian, choose any
$\epsilon > 0$. Since Lemma 1 states that
 $\Gamma^-$ is uniformly integrable on $SV$, we can choose
 $\delta > 0$ so that \[ A \subseteq SV \ \mathrm{and} \ \mu\left(A\right) < \delta \Longrightarrow  \forall f \in \Gamma, \ \ \
 \left\Vert f^-\right\Vert _{1; A} < \frac{\epsilon}{4  \left\Vert j\right\Vert _{\infty ; A}}. \]
  By Lusin's theorem,
  we may choose a compact $F  \subseteq S V$
  where \[\mu\left(S V\right) - \mu\left(F\right) < \frac{\delta}{2}\] and $j |_F$ is continuous.
  Moreover, by the Tietze extension theorem, we may extend $j$ to a continuous function $g$ on $S V$,
   which provides us with a continuous $g : S V \rightarrow \mathbb{C}$
 and a compact $F \subseteq S V$ such that
  $j = g$ on $ F$, where \[ \mu\left(S V\right) - \mu\left(F\right) < \frac{\delta}{2}.\]

\indent Since $g$ has compact support,  for any $x \in G$, there
exists a relative neighborhood $N_x \subseteq S$ of $x$ such that
\[x' \in N_x \Longrightarrow |j\left(xy\right) - j\left(x'y\right)| = |g\left(xy\right) - g\left(x' y\right)| <
\frac{\epsilon}{2M}\] for all $y \in \left(x^{-1} F \cap {x'}^{-1}
F\right)$.

 \indent Thus fixing $x$ and fixing any such $x'$, we find that

\begin{eqnarray}\left\Vert \left(j\ast f\right)\left(x\right) -
\left(j\ast f\right)\left(x'\right)\right\Vert _B  &\leq&  \int_V
|j\left(xy\right) - j\left(x' y\right)| \left\Vert
f\left(y^{-1}\right) \right\Vert _B \, dy \nonumber
\\ &\leq& \int_{S V \backslash \left(x^{-1} F \cap {x'}^{-1}
F\right)}
 |j\left(xy\right) - j\left(x' y\right)| \left\Vert f\left(y^{-1}\right) \right\Vert _B \, dy  \nonumber\\ &+& \int_{\left(x^{-1} F \cap {x'}^{-1} F\right) }
 |j\left(xy\right) - j\left(x' y\right)| \left\Vert f\left(y^{-1}\right) \right\Vert _B \, dy \nonumber\\ &<& \frac{\epsilon}{2} +
\int_{\left(x^{-1} F \cap {x'}^{-1} F\right) }
 |j\left(xy\right) - j\left(x' y\right)| \left\Vert f\left(y^{-1}\right) \right\Vert _B \, dy \nonumber\\ &<& \frac{\epsilon}{2} +
 \frac{\epsilon}{2M} \int_{\left(x^{-1} F \cap {x'}^{-1} F\right) }
  \left\Vert f\left(y^{-1}\right)\right\Vert _B \, dy  < \epsilon.\nonumber \end{eqnarray}

 \indent Therefore, we have \[ \forall f \in \Gamma, \ \ \ \left\Vert \left(j\ast f\right)\left(x\right) -
\left(j\ast f\right)\left(x'\right)\right\Vert _B < \epsilon. \]
Since $x \in S$ was arbitrary, it follows  that $j\ast\Gamma$ is
equi-continuous on
$S$.\\

\indent Now we check that the set
$\left(j\ast\Gamma\right)\left(s\right) = \{\left(j \ast
f\right)\left(s\right) : f \in \Gamma \}$ is relatively norm compact
in $B$ for every $s \in S$. Although Dinculeanu proved this for $G$
abelian, and although the proof is nearly identical for $G$
non-abelian, we nevertheless present the proof for the sake of
completeness. Assume that $\mu\left(A\right) < \infty
\Longrightarrow \{\int_A f^- \left(y\right) \, dy : f \in \Gamma\}$
is relatively norm compact in $B$ and that $G$ is not necessarily
abelian.

 Let  $j$ be
a simple function of the form $j\left(y\right) = \sum {\alpha}_i
\chi_{A_i}\left(y\right)$ where each $A_i$ has finite measure and
where the sum is taken over finitely many $i$. Then we have
\[\left(j\ast f\right)\left(x\right) = \sum {\alpha}_i \int_G
\chi_{A_i}\left(xy\right) f\left(y^{-1}\right) \, dy = \sum
{\alpha}_i \int_{x^{-1} A_i } f\left(y^{-1}\right) \, dy.\] Thus,
since each $\{\int_{x^{-1} A_i } f\left(y^{-1}\right) \, dy : f \in
\Gamma\}$ is relatively norm compact in $B$, we have that
$\left(j\ast\Gamma\right)\left(x\right) = \{\sum {\alpha}_i
\int_{x^{-1} A_i } f\left(y^{-1}\right) \, dy : f \in \Gamma\}$ is
relatively norm compact in $B$.

\indent If $j$ is not a simple function, then let $j_n \rightarrow
j$ in the $L^{p'}$ norm where $j_n$ is a simple function.  We show
that $\left(j_n \ast f\right)\left(x\right) \rightarrow \left(j \ast
f\right)\left(x\right)$ uniformly for  $f \in \Gamma.$   To that
end, again let $ 0 < m < \triangle\left(y^{-1}\right) $ on $V$, so
from
 H\"{o}lder's inequality we have
\begin{eqnarray} \left\Vert \left(j_n \ast f\right)\left(x\right) -
\left(j \ast f\right)\left(x\right)\right\Vert _B & \leq& \int_G
|j_n\left(y\right) - j\left(y\right)| \left\Vert
f\left(y^{-1}x\right)\right\Vert \, dy \nonumber\\& \leq &\left\Vert
j_n - j\right\Vert _{p' ; V} \left(\int_V \left\Vert f\left(y^{-1}
x\right)\right\Vert _B ^p \, dy\right)^{\frac{1}{p}} \nonumber\\& <&
\frac{\left(\triangle\left(x\right)\right)^{\frac{1}{p}}}{m}
\left\Vert j_n - j\right\Vert _{p' ; V} \left(\int_V \left\Vert
f\left(y\right)\right\Vert _B ^p \, dy \right)^{\frac{1}{p}}
\nonumber\\& \leq& \frac{M}{m} \left({\triangle\left(x\right)}
\right)^{\frac{1}{p}} \left\Vert j_n - j\right\Vert _{p' ;
V}\nonumber
\end{eqnarray}

where $\Gamma$ is norm bounded by $M$.

 \indent  Thus,
$\left(j_n \ast f\right)\left(x\right) \rightarrow \left(j \ast
f\right)\left(x\right)$ uniformly for $f \in \Gamma.$ Note that if
$G$ is abelian, then unimodularity implies that $j$ does not need
compact support. Also note that if $j$ has compact support $V$, then
we may assume that $A \subseteq V$ so that $A$ is relatively
compact, and therefore each $\int_A f^- \left(y\right) \, dy $
exists.

Finally, by Lemma $4$, each $\left(j\ast\Gamma\right)\left(x\right)$
for $x \in K$ is norm compact in $B$.  Therefore, by Lemma $6$, we
are done.

\end{pf}\qed

We note that if $G$ is abelian or compact, then unimodularity
implies that the following two statements are equivalent:

1) $A \subseteq G \ \mathrm{and} \ \mu\left(A\right) < \infty
\Longrightarrow \{\int_A f^- \left(y\right) \, dy : f \in \Gamma \}$
is relatively norm compact in $B$.

2) $A \subseteq G \ \mathrm{and} \ \mu\left(A\right) < \infty
\Longrightarrow \{\int_A f \left(y\right) \, dy : f \in \Gamma\}$ is
relatively norm compact in $B$.

\section{Main Result}

Now we may state our principal theorem.\\

\begin{thm}
Let $G$ be a locally compact  group with left invariant Haar measure
$\mu$ and $B$ be a Banach space.
  A bounded subset $\Gamma$ of the Bochner space $L^p\left(G ; B\right)$ is relatively norm compact
  iff

1) for any $A \subseteq G $ with $\ \mu \left(A\right) < \infty$,
the set $\{\int_A f^- \left(y\right) \, dy : f \in \Gamma\}
\subseteq B$ is relatively
norm compact in $B$,\\

 and one of the following is true:

2a) for any $\epsilon > 0$, there exists a compact $S \subseteq G$
and $j \in
 L^{p'}\left(G\right)$ such that \[ \forall f \in \Gamma, \ \ \ \left\Vert f\right\Vert
_{p ; G \backslash S} < \epsilon\] and \[ \forall f \in \Gamma, \ \
\ \left\Vert f - j\ast f\right\Vert _{p ; S} < \epsilon\] where $j$
has compact support if either $G$ is nonabelian or $p =1$,

2b) for any $\epsilon > 0$, there exists a compact $S \subseteq G$
and a neighborhood $\Theta$ of the identity such that \[ \forall f
\in \Gamma, \ \ \ \left\Vert f\right\Vert _{p ; G \backslash S} <
\epsilon\] and
\[\forall f \in \Gamma, \ \ \ h \in \Theta \Longrightarrow \left\Vert T^h f - f\right\Vert _{p ; S} <
\epsilon\] where  $T^h$ is the translation operator defined by
$\left(T^h f\right)
\left(x\right) = f\left(h^{-1} x\right)$.\\\\

\end{thm}

\begin{pf} First we prove necessity for the general case of a locally compact
group and $1 \leq p < \infty$. We show that both 2a) and 2b) hold.
To that end, let $\epsilon > 0$. Since $\Gamma$ is totally bounded,
there exists $\{f_1,\ldots,f_N\} \subseteq \Gamma$ such that
\[\Gamma \subseteq \bigcup_{i=1}^N B\left(f_i, \frac{\epsilon}{3}\right)\]
where $B\left(f_i, \frac{\epsilon}{3}\right) $ denotes a ball of
radius $\frac{\epsilon}{3}$ at $f_i$ in $L^p\left(G, B\right).$

\indent Choose a compact $S \subseteq G$ such that  \[\left\Vert
f_i\right\Vert _{p ; G \backslash S} < \frac{\epsilon}{2}\] for all
$i \in \{1,\ldots,N\}$.
 Then, we have
 \[ \forall f \in \Gamma, \ \ \ \left\Vert f\right\Vert _{p ; G \backslash S} \leq \left\Vert f - f_k\right\Vert _{p ; G \backslash
 S} + \left\Vert f_k\right\Vert _{p ; G \backslash S} < \epsilon, \]
 where $k \in \{1,\ldots,N\}$ corresponds to $f$.

\indent For property 2a), let $\mathfrak{V}$ be a basis of
relatively compact neighborhoods of the identity of $G$. For each $V
\in \mathfrak{V}$, let $u_{V}=\mu (V)^{-1}\chi _{V}$.
 We show that $\lim_{V \in
\mathfrak{V} } \left\Vert u_V \ast f - f\right\Vert _{p ; G} = 0$
uniformly in $\Gamma$.  This holds for any $f \in
\{f_1,\ldots,f_N\}$,
 so choose $V'$ such that \[
\forall f \in \{f_1,\ldots,f_N\}, \ \ \ V \subseteq V'
\Longrightarrow \left\Vert u_V \ast f - f\right\Vert _{p ;  G} <
\frac{\epsilon}{3}.\] \indent

Since $\left\Vert u_V\right\Vert _{1;G} = 1$, a standard computation
using H\"{o}lder's inequality and left invariance shows that
$\left\Vert u_V \ast f - u_V \ast f_k \right\Vert _{p ; G} \leq
\left\Vert f - f_k\right\Vert _{p ; G}$. Therefore, we have
\[
\begin{array}{lll}
\forall f\in \{f_{1},...,f_{N}\},\ \ \ V\subseteq V^{\prime } &
\Longrightarrow  &
\left\Vert u_{V}\ast f-f\right\Vert _{p;G} \\
& \leq  & \left\Vert u_{V}\ast f-u_{V}\ast f_{k}\right\Vert
_{p;G}+\left\Vert u_{V}\ast f_{k}-f_{k}\right\Vert _{p;G} \\
&  & \qquad +\left\Vert f_{k}-f\right\Vert _{p;G} \\
& \leq  & 2\left\Vert f-f_{k}\right\Vert _{p;G}+\left\Vert
f_{k}-j\ast
f_{k}\right\Vert _{p;G} \\
& < & \epsilon ,%
\end{array}%
\]where $k \in \{1,\ldots,N\}$
corresponds to $f$.

\indent For property 2b), choose a neighborhood $\Theta$ of the
identity with the property
\[\forall k \in \{1,\ldots,N\},\ \ \ y \in \Theta \Longrightarrow \int_{G} \left\Vert f_k\left(x\right) -
  f_k\left(y^{-1} x\right)\right\Vert _B ^p \, dx <
  \left({\frac{\epsilon}{3}}\right)^p.\]
Thus, we have \[
\begin{array}{lll}
\forall f\in \Gamma ,\ \ \ y\in \Theta  & \Longrightarrow  & \left(
\int_{G}\left\Vert f(x)-f(y^{-1}x)\right\Vert _{B}^{p}dx\right) ^{1/p} \\
& \leq  & \left( \int_{G}\left\Vert f(x)-f_{k}(x)\right\Vert
_{B}^{p}dx\right) ^{1/p}+\left( \int_{G}\left\Vert
f_{k}(x)-f_{k}(y^{-1}x)\right\Vert _{B}^{p}dx\right) ^{1/p} \\
&  & \qquad +\left( \int_{G}\left\Vert
f_{k}(x)-f(y^{-1}x)\right\Vert
_{B}^{p}dx\right) ^{1/p} \\
& < & \epsilon
\end{array}%
\]
by left invariance.

\indent Finally, if $\mu\left(A\right) < \infty$ and $G$ is abelian,
then the continuity of the map $f \mapsto \int_A f^-\left(x\right)
\, dx$ of $L^p\left(G ; B\right)$ into $B$ implies that $\{\int_A
f^-\left(x\right) \, dx : f \in \Gamma\} \subseteq B$ is norm
compact in $B$.  If $G$ is nonabelian, then as we noted in the proof
of Lemma 7, we only need to consider the case where $A$ is
relatively compact. Therefore, the continuity of the map $f \mapsto
\int_A f^-\left(x\right) \, dx$ easily follows.

Now we prove sufficiency. Assume that properties 1) and 2a) hold.
Let $\{f_n\} \subseteq \Gamma$ and choose any $\epsilon > 0$. Choose
a compact $S \subseteq G$ and $j \in L^{p'} \left(G\right)$  so that
\[\forall f_n \in \Gamma, \ \ \ \left\Vert f_n\right\Vert _{p ;  G
\backslash S} < \frac{\epsilon}{6}\]and
\[\forall f_n \in \Gamma, \ \ \ \left\Vert f_n - j \ast
{f_n}\right\Vert _{p ; S} < \frac{\epsilon}{3}.
\] Morever, since $G$ is locally compact, we may assume that $S$ contains the identity. Therefore,
 Lemma $7$ implies that $\left(j \ast \Gamma \right)| _S $ is
relatively compact in $C_B\left(S\right)$, and so passing to a
subsequence if necessary (and assuming without loss of generality
that each $j \ast {f_n}$ is restricted to $S$), we have that $j \ast
{f_n} \rightarrow f$ in
 $C_B\left(S\right)$ for some $f \in C_B\left(S\right)$. Define
 $f$ to be zero on $G \backslash S$.

\indent By our hypothesis and the fact that $ j \ast {f_n}
\rightarrow f$ in
 $C_B\left(S\right)$, we have that for all large enough $n \in \mathbb{N}$
,
\[
\begin{array}{lll}
\forall f\in \Gamma ,\ \ \ \left\Vert f_{n}-f\right\Vert _{p;G} &
\leq & \left\Vert f_{n}-f\right\Vert _{p;G\setminus S}+\left\Vert
f_{n}-f\right\Vert _{p;S} \\
& \leq  & \epsilon /3+\left\Vert f_{n}-f\right\Vert _{p;S} \\
& \leq  & \epsilon /3+\left\Vert f_{n}-j\ast f_{n}\right\Vert
_{p;S}+\left\Vert j\ast f_{n}-f\right\Vert _{p;S} \\
& \leq  & \epsilon
\end{array}%
\]
 Hence, $\Gamma$ is relatively norm compact.

\indent Finally, we show that property 2b) implies property 2a).
 Choose a neighborhood $\Theta'$ of the identity such that \[ y \in \Theta'
 \Longrightarrow \int_{S}
\left\Vert f\left(x\right) - f\left(y^{-1} x\right) \right\Vert _B
^p \, dx < {\epsilon}^p .\] Since $G$ is locally compact, we may
assume $\Theta'$ is relatively compact, and so $\Theta'$ has nonzero
finite measure. Let $j = \mu\left(\Theta'\right)^{-1} \chi
_{\Theta'}$.

\indent Then for
 any $f \in \Gamma$, we find, using H\"{o}lder's inequality and Fubini's theorem, that
\begin{eqnarray} \forall f \in \Gamma, \ \ \ \left\Vert f -
j \ast f\right\Vert _{p ; S}^p &= & \int_S \left\Vert
f\left(x\right) - \int_{ \Theta'} j\left(y\right) f\left(y^{-1}
x\right) \, dy \right\Vert _B ^p \, dx  \nonumber
\\& \leq &\int_S \left(\int_{\Theta'} j\left(y\right)^\frac{1}{p'}
j\left(y\right)^\frac{1}{p} \left\Vert f\left(x\right) -
f\left(y^{-1} x\right) \right\Vert _B \, dy \right)^p \, dx
\nonumber
\\& \leq &\int_S \left(\int_{\Theta'} j\left(y\right) \left\Vert
f\left(x\right) - f\left(y^{-1} x\right) \right\Vert _B ^p \, dy
\right)dx \nonumber  \\& =& \int_{\Theta'} \left(\int_{S} \left\Vert
f\left(x\right) - f\left(y^{-1} x\right) \right\Vert _B ^p \, dx
\right)j\left(y\right) \, dy .\nonumber
\end{eqnarray}

\indent But for all $y \in \Theta'$, we have  $\int_S \left\Vert
f\left(x\right) - f\left(y^{-1} x\right) \right\Vert _B ^p \, dx <
{\epsilon}^p $, and so $\left\Vert f - j \ast f\right\Vert _{p ; S}
< \epsilon$.
\end{pf}\qed

\section{Final Comments}

 \indent We make two final comments. First, Theorem $1$ was proven
  in ~\cite{Book7} for the case where $B = \mathbb{C}$, $G$ is an arbitrary
 locally compact group,  and each $u_V = \mu\left(V\right)^{-1} \chi_V $ where $V$ is a nonempty neighborhood
  of the identity, so that our
  theorem generalizes this result in two different directions.

  \indent Also, in the case $B = \mathbb{C}$ and $G = \mathbb{R}^n$,  Shurenkova and
  Buldygin ~\cite{Book1}
  proved Theorem $1$  for the case when $j$ is a function
  of the form $\omega_h$ for some $h > 0$, where
   $\{\omega_h : h \in \mathbb{R}^+\}$ is any family of real functions
    in $ \subseteq L^1\left(\mathbb{R}^n\right)
 \cap L^{p'}\left(\mathbb{R}^n\right)$ satisfying the following:

 1) $h > 0, \left\Vert \omega_h\right\Vert _{1; \mathbb{R}^n} = 1,\text{and } \omega_h $ is continuous
 a.e. on $\mathbb{R}^n$,

 2) for any $\delta > 0$, $\lim_{h \rightarrow 0} \left\Vert \omega_h\right\Vert _{p ;  \{|x| > \delta\}} = 0$,

 3) $\left\Vert \omega_h\right\Vert _{1;\mathbb{R}^n} = 1$

(here of course, we do not consider the property 3b) in the
statement of Theorem 1.)

Thus, considering the well-known fact that
$L^p\left(\mathbb{R}^n\right)$ functions can be arbitrarily
approximated (in the $L^p$ norm) by such functions ~\cite[chap. 9,
Sect. 2]{Book9}, it is clear that our theorem is a significant
generalization of Shurenkova and Buldygin's results.

\begin{ack}The author
would like to thank the entire NJIT mathematics department, and
especially Professor Denis Blackmore, for their time and very
generous help.
\end{ack}


\begin{thebibliography}{8}


\bibitem{Book1}
{\sc Buldygin, V.} and {\sc Shurenkova, A.} \emph{On a Criterion of
Compactness in the Space $L_p$}, Theory of Probability and
Mathematical Statistics, No. 52, 1996, 33 - 37.

\bibitem{Book2}
{\sc Dinculeanu, N.} \emph{On Kolmogorov-Tamarkin and M. Riesz
Compactness Criteria in Function Spaces Over a Locally Compact
Group}, Journal of Mathematical Analysis and Applications, No. 89,
1982, 67 - 85.

\bibitem{Book3}
{\sc Dinculeanu, N.} \emph{Strong Additivity, Absolute Continuity,
and Compactness in Spaces of Measures}, Journal of Mathematical
Analysis and Applications, No. 45, 1974, 172 - 188.


\bibitem{Book4}
{\sc Hille, E.} and {\sc Phillips, R.}\emph{Functional Analysis and
Semi-Groups}, American Mathematical Society., 1957


\bibitem{Book5}
{\sc Loomis, L.} \emph{An Introduction to Abstract Harmonic
Analysis}, Van Nostrand Company, inc., 1953.


\bibitem{Book6}
{\sc Munkres, J.} \emph{Topology} Prentice Hall, 1999.



\bibitem{Book7}
{\sc Teleman, S.} \emph{Sur Les Ensembles Compacts de Fonctions
Sommables}, Rev. Math. Pures Appl. No. 6, 1961, 659?684.

\bibitem{Book8}
{\sc Welland, R.} and {\sc Goes, S.} \emph{Compactness Criteria for
K\"{o}the Spaces} Math Ann. No. 188, 1970, 251 - 269.

\bibitem{Book9}
{\sc Wheeden, R.} and {\sc Zygmund, A.} \emph{Measure and Integral},
Marcel Dekker Inc., 1977.

\end{thebibliography}
\end{document}